\documentclass[12pt]{article}

\usepackage{amsmath,amssymb}
\usepackage[applemac]{inputenc}
\usepackage{amsmath,amssymb,fullpage}
\usepackage{showlabels}
\usepackage{color}

\newenvironment{myenumerate}{%

\begin{enumerate}}{\end{enumerate}}

\newcommand{\dproof}{\noindent {Proof.} \quad}
\newcommand{\fproof}{\hfill $\square$ \bigskip}

\newtheorem{definition}{Definition}[section]

\newtheorem{theorem}[definition]{Theorem}
\newtheorem{problem}[definition]{Problem}

\newtheorem{coro}[definition]{Corollary}

\numberwithin{equation}{section}

\def\RB{\mathbb{R}}

\def\1B{\text{1\!\!I}}

\begin{document}
\date{14 August 2016 }
\title{Optimal insider control of stochastic partial differential equations}

\author{
Olfa Draouil$^{1}$ and Bernt \O ksendal$^{2,3}$}

\footnotetext[1]{Department of Mathematics, University of Tunis El Manar, Tunis, Tunisia.\\
Email: {\tt olfadraouil@hotmail.fr}}

\footnotetext[2]{Department of Mathematics, University of Oslo, P.O. Box 1053 Blindern, N--0316 Oslo, Norway.\\
Email: {\tt oksendal@math.uio.no}}

\footnotetext[3]{This research was carried out with support of the Norwegian Research Council, within the research project Challenges in Stochastic Control, Information and Applications (STOCONINF), project number 250768/F20.}

\maketitle



\begin{abstract} We study the problem of optimal inside control of an SPDE (a stochastic evolution equation) driven by a Brownian motion and a Poisson random measure.  Our optimal control problem is new in two ways:
\begin{itemize}
\item
(i) The controller has access to inside information, i.e. access to information about a future state of the system,
\item
(ii) The integro-differential operator of the SPDE might depend on the control.
\end{itemize}

In the first part of the paper, we formulate a sufficient and a necessary maximum principle for this type of control problem, in two cases:
\begin{itemize}
\item
The control is allowed to depend both on time $t$ and on the space variable $x$.
\item
The control is not allowed to depend on $x$.
\end{itemize}

In the second part of the paper, we apply the results above to the problem of optimal control of an SDE system when the inside controller has only noisy observations of the state of the system. Using results from nonlinear filtering, we transform this noisy observation SDE inside control problem into a \emph{full observation} SPDE insider control problem.\\
The results are illustrated by explicit examples.

\end{abstract}

\section{Introduction}
In this paper we consider an optimal control problem for a stochastic process $Y(t,x)=Y^{u,Z}(t,x)=Y(t,x,Z)=Y(t,x,z)|_{z=Z}$ defined as the solution of a stochastic partial differential equation (SPDE) given by
\begin{align}\label{eq1.3}
&dY(t,x)=[A_{u(t,x,Z)}Y(t,x)+a(t,x,Y(t,x),u(t,x,Z),Z)]dt+b(t,x,Y(t,x),u(t,x,Z),Z)dB(t)\nonumber\\
&+\int_{\mathbb{R}} c(t,x,Y(t,x),u(t,x,Z),Z,\zeta)\tilde{N}(dt,d\zeta); \quad (t,x)\in (0, T)\times D.
\end{align}
The boundary conditions are
\begin{equation}
Y(0,x)=\xi(x),\quad x\in D
\end{equation}
\begin{equation}
Y(t,x)=\theta(t,x); \quad (t,x)\in [0,T]\times \partial D.
\end{equation}
Here $B(t)$ and $\tilde{N}(dt,d\zeta)$ is a Brownian motion and an independent compensated Poisson random measure, respectively, jointly defined on a filtered probability space $(\Omega, \mathbb{F}=\{ \mathcal{F}_t \}_{t \geq 0},\mathbf{P})$ satisfying the usual conditions. $T>0$ is a given constant, $D \subset \mathbb{R}$ is a given open set, and $\partial D$ denotes the boundary of $D$.
The process $u(t,x)=u(t,x,z)_{z=Z}$ is our insider control process, where $Z$ is a given $\mathcal{F}_{T_0}$-measurable random variable for some $T_0>0$ , representing the inside information available to the controller.
\\

The operator $A_{u}$ is a linear integro-differential operator acting on $x$, with parameter $u$, and the expression
$A_{u(t,x,Z)}Y(t,x))$ means $A_u Y(t,x,Z)|_{u=u(t,x,Z)}$.\\

We interpret the equation \eqref{eq1.3} for $Y$ in the weak sense. By this we mean that $Y(t,\cdot)$ satisfies the equation
\begin{align}
&(Y(t,\cdot),\phi)_{L^2(D)} = (\xi,\phi)_{L^2(D)}+ \int_0^t(Y(s,\cdot),A_u^*\phi)_{L^2(D)}ds+\int_0^t(a(s,Y(s,\cdot),\cdot),\phi)_{L^2(D)}ds \nonumber\\
&+ \int_0^t (b(s,Y(s,\cdot),\cdot),\phi)_{L^2(D)}dB(s)+ \int_0^t\int_{\mathbb{R}}c(s,Y(s,\cdot),\zeta,\cdot),\phi)_{L^2(D)}\tilde{N}(ds,d\zeta),
\end{align}
for all smooth functions $\phi$ with compact support in $D$. Here
\begin{equation}
(\psi,\phi)_{L^2(D)}=\int_{D}\psi(x)\phi(x)dx
\end{equation}
is the $L^2$ inner product on $D$ and $A_u^*$ is the adjoint of the operator $A_u$, in the sense that
\begin{equation}\label{adjoint}
(A_u\psi,\phi)_{L^2(D)} = (\psi,A_u^* \phi)_{L^2(D)}
\end{equation}
for all smooth $L^2$ functions $\psi,\phi$ with compact support in D. It can be proved that the It\^ o formula can be applied to such SPDEs. See \cite{Par}, \cite{PR}.

 We assume that the inside information is of \emph{initial enlargement} type. Specifically, we assume  that the inside filtration $\mathbb{H}$ has the form

\begin{equation}\label{eq1.1}
 \mathbb{H}= \{ \mathcal{H}_t\}_{0\leq t \leq T}, \text{ where } \mathcal{H}_t = \mathcal{F}_t \vee \sigma(Z)
\end{equation}
for all $t$, where $Z$ is a given $\mathcal{F}_{T_0}$-measurable random variable, for some $T_0 > 0$ (constant).
Here and in the following we use the right-continuous version of  $ \mathbb{H}$, i.e. we put
$\mathcal{H}_{t}= \mathcal{H}_{t^+}=\bigcap_{s>t}\mathcal{H}_s.$

We also assume that the \emph{Donsker delta functional} of $Z$ exists (see below). This assumption implies that the Jacod condition holds, and hence that $B(\cdot)$ and $N(\cdot,\cdot)$ are semimartingales with respect to $\mathbb{H}$. See e.g. \cite{DO3} for details.
We assume that the value at time $t$ of our insider control process $u(t,x)$ is allowed to depend on both $Z$ and $\mathcal{F}_t$. In other words, $u(.,x)$ is assumed to be $\mathbb{H}$-adapted. Therefore it has the form
\begin{equation}\label{eq1.2}
    u(t,x,\omega) = u_1(t, x,Z, \omega)
\end{equation}
for some function $u_1 : [0, T]\times D \times \mathbb{R}\times\Omega\rightarrow \mathbb{R}$
such that $u_1(. ,x, z)$ is $\mathbb{F}$-adapted for each $(x,z) \in D \times \mathbb{R}$. For simplicity (albeit with some abuse of notation) we will in the following write $u$ instead of $u_1$.

Let $\mathbb{U}$ denote the set of admissible control values.We assume that the functions
\begin{align}
a(t,x,y,u,z)&=a(t,x,y,u,z,\omega):
[0,T]\times D \times \mathbb{R} \times \mathbb{U}\times \mathbb{R} \times \Omega \mapsto \mathbb{R}\nonumber\\
b(t,x,y,u,z)&=b(t,x,y,u,z,\omega):
[0,T]\times D\times \mathbb{R} \times \mathbb{U} \times \mathbb{R} \times \Omega \mapsto \mathbb{R}\nonumber\\
c(t,x,y,u,z,\zeta)&=c(t,x,y,u,z,\zeta,\omega):
[0,T]\times D \times \mathbb{R} \times \mathbb{U} \times \mathbb{R} \times \mathbb{R} \times \Omega \mapsto \mathbb{R}\nonumber\\
\end{align}
are given bounded $C^1$ functions with respect to $y$ and $u$ and adapted processes in $(t,\omega)$ for each given $x,y,u,z,\zeta$.
Let $\mathcal{A}$ be a given family of admissible $\mathbb{H}-$adapted controls $u$.
The \emph{performance functional} $J(u)$ of a control process $u \in \mathcal{A}$ is defined by
\begin{equation}\label{eq1.4}
J(u)= \mathbb{E}[\int_0^T (\int_Dh(t,x,Y(t,x),u(t,x,Z),Z)dx)dt +\int_Dk(x,Y(T,x),Z)dx],
\end{equation}
where \begin{align}
&h(t,x,y,u,z): [0,T] \times D \times\mathbb{R}\times \mathbb{U} \times\mathbb{R} \mapsto \mathbb{R}\nonumber\\
&k(x,y,z):D\times\mathbb{R} \times \mathbb{R} \mapsto \mathbb{R}
\end{align}
are given bounded functions, $C^1$ with respect to $y$ and $u$. The functions $h$ and $k$ are called the \emph{profit rate density} and \emph{terminal payoff density}, respectively. For completeness of the presentation we allow these functions to depend explicitly on the future value $Z$ also, although this would not be the typical case in applications. But it could be that $h$ and $k$ are influenced by the future value $Z$ directly through the action of an insider, in addition to being influenced indirectly through the control process $u$ and the corresponding state process $Y$.\\

\begin{problem}
Find $u^{\star} \in\mathcal{A}$ such that
\begin{equation}\label{eq1.5}
    \sup_{u\in\mathcal{A}}J(u)=J(u^{\star}).
\end{equation}
\end{problem}

\section{The Donsker delta functional}
To study this problem we adapt the technique of the paper \cite{DO1} to the SPDE situation and we combine this with the method for  optimal control of SPDE developed in \cite{O1}, \cite{OPZ} and \cite{OS1}. We first recall briefly the definition and basic properties of the Donsker delta functional:

\begin{definition}
Let $Z:\Omega\rightarrow\mathbb{R}$ be a random variable which also belongs to $(\mathcal{S})^{\ast}$. Then a continuous functional
\begin{equation}\label{donsker}
    \delta_{Z}(.): \mathbb{R}\rightarrow (\mathcal{S})^{\ast}
\end{equation}
is called a Donsker delta functional of $Z$ if it has the property that
\begin{equation}\label{donsker property }
    \int_{\mathbb{R}}g(z)\delta_{Z}(z)dz= g(Z) \quad a.s.
\end{equation}
for all (measurable) $g : \mathbb{R} \rightarrow \mathbb{R}$ such that the integral converges.
\end{definition}

For example, consider the special case when $Z$ is a first order chaos random variable of the form
\begin{equation}\label{eq2.5}
    Z = Z (T_0); \text{ where } Z (t) =\int_0^t\beta(s)dB(s)+\int_0^t\int_{\mathbb{R}}\psi(s,\zeta)\tilde{N}(ds,d\zeta), \mbox{ for } t\in [0,T_0]
\end{equation}
for some deterministic functions $\beta \neq 0, \psi$ such that
\begin{equation}\label{}
    \int_0^{T_0} \{ \beta^2(t)+\int_{\mathbb{R}}\psi^2(t,\zeta)\nu(d\zeta)\} dt<\infty \text{ a.s. }
\end{equation}
and for every $\epsilon >0$ there exists $\rho > 0$ such that
\begin{equation}\label{eq8.4}
\int_{\mathbb{R} \setminus (-\epsilon,\epsilon)} e^{\rho  \zeta} \nu(d\zeta) < \infty.\nonumber\\
\end{equation}

This condition implies that the polynomials are dense in $L^2(\mu)$, where $d\mu(\zeta)=\zeta^2 d\nu(\zeta)$. It also guarantees that the measure $\nu$ integrates all polynomials of degree $\geq 2$.\\
In this case it is well known (see e.g. \cite{MOP}, \cite{DiO1}, Theorem 3.5, and \cite{DOP},\cite{DiO2}) that the Donsker delta functional exists in $(\mathcal{S})^{\ast}$ and is given
by
\begin{eqnarray}\label{eq2.7}
   \delta_Z(z)&=&\frac{1}{2\pi}\int_{\mathbb{R}}\exp^{\diamond}\big[ \int_0^{T_0}\int_{\mathbb{R}}(e^{ix\psi(s,\zeta)}-1)\tilde{N}(ds,d\zeta)+ \int_0^{T_0}ix\beta(s)dB(s)  \nonumber\\
   &+&  \int_0^{T_0}\{\int_{\mathbb{R}}(e^{ix\psi(s,\zeta)}-1-ix\psi(s,\zeta))\nu(d\zeta)-\frac{1}{2}x^2\beta^2(s)\}ds-ixz\big]dx,
\end{eqnarray}
where $\exp^{\diamond}$ denotes the Wick exponential.
Moreover, we have for $t< T_0$
\begin{align}
&\mathbb{E}[\delta_Z(z)|\mathcal{F}_t]\nonumber\\
=&  \frac{1}{2\pi}\int_{\mathbb{R}}\exp\big[\int_0^t\int_{\mathbb{R}}ix\psi(s,\zeta)\tilde{N}(ds,d\zeta) +\int_0^t ix\beta(s)dB(s)\\
&+\int_t^{T_0}\int_{\mathbb{R}}(e^{ix\psi(s,\zeta)}-1-ix\psi(s,\zeta))\nu(d\zeta)ds-\int_t^{T_0}\frac{1}{2}x^2\beta^2(s)ds-ixz\big]dx.
\end{align}

If $D_t$ and $D_{t,\zeta}$ denotes the \emph{Hida-Malliavin derivative} at $t$ and $t,\zeta$  with respect to $B$ and $\tilde{N}$, respectively, we have
\begin{align}
&\mathbb{E}[D_t\delta_Z(z)|\mathcal{F}_t]=\nonumber\\
&\frac{1}{2\pi}\int_{\mathbb{R}}\exp\big[\int_0^t\int_{\mathbb{R}}ix\psi(s,\zeta)\tilde{N}(ds,d\zeta) +\int_0^t ix\beta(s)dB(s)\nonumber\\
&+\int_t^{T_0}\int_{\mathbb{R}}(e^{ix\psi(s,\zeta)}-1-ix\psi(s,\zeta))\nu(d\zeta)ds-\int_t^{T_0}\frac{1}{2}x^2\beta^2(s)ds-ixz\big]ix\beta(t)dx
\end{align}
   and
\begin{align}
&\mathbb{E}[D_{t,z}\delta_Z(z)|\mathcal{F}_t]=\nonumber\\
& \frac{1}{2\pi}\int_{\mathbb{R}}\exp\big[\int_0^t\int_{\mathbb{R}}ix\psi(s,\zeta)\tilde{N}(ds,d\zeta) +\int_0^t ix\beta(s)dB(s)\nonumber\\
&+\int_t^{T_0}\int_{\mathbb{R}}(e^{ix\psi(s,\zeta)}-1-ix\psi(s,\zeta))\nu(d\zeta)ds-\int_t^{T_0}\frac{1}{2}x^2\beta^2(s)ds-ixz\big](e^{ix\psi(t,z)}-1)dx.
\end{align}

For more information about the Donsker delta functional, Hida-Malliavin calculus and their properties, see \cite{DO1}.

From now on we assume that $Z$ is a given random variable which also belongs to $(\mathcal{S})^{\ast}$, with a Donsker delta functional $\delta_Z(z)\in(\mathcal{S})^{\ast}$
satisfying
\begin{equation}
\mathbb{E}[\delta_Z(z)|\mathcal{F}_T] \in \mathbf{L}^2(\mathcal{F}_T,P)
\end{equation}
and
\begin{equation}
\mathbb{E}[\int_0^T(\mathbb{E}[D_t\delta_Z(z)|\mathcal{F}_t])^2 dt ]<\infty, \text{ for all } z.
\end{equation}

\section{Transforming the insider control problem to a related parametrized non-insider problem}
Since $Y(t,x)$ is $\mathbb{H}$-adapted, we get by using the definition of the Donsker delta functional $\delta_Z(z)$ of $Z$ that
\begin{equation}\label{eq1.6}
Y(t,x)=Y(t,x,Z)=Y(t,x,z)_{z=Z}=\int_{\mathbb{R}}Y(t,x,z)\delta_Z(z)dz
\end{equation}
for some $z$-parametrized process $Y(t,x,z)$ which is $\mathbb{F}$-adapted for each $x,z$.
Then, again by the definition of the Donsker delta functional we can write, with
$A_u=A_{u(s,x,Z)}=A_{u(s,x,z)_{z=Z}}$,
\begin{align}\label{eq1.7}
&Y(t,x)= \xi(x,Z) +\int_0^t[A_{u}Y(s,x)+ a(s,x,Y(s,x),u(s,x,Z),Z)]ds + \int_0^t b(s,x,Y(s,x),u(s,x,Z),Z)dB(s)\nonumber\\
&+\int_0^t \int_{\mathbb{R}} c(s,x,Y(s,x),u(s,x,Z),Z,\zeta)\tilde{N}(ds,d\zeta)\nonumber\\
&=\xi(x,z)_{z=Z}+\int_0^t[A_uY(s,x,z)+ a(s,x,Y(s,x,z),u(s,x,z),z)]_{z=Z}ds \nonumber\\
&+ \int_0^t b(s,x,Y(s,x,z),u(s,x,z),z)_{z=Z}dB(s)\nonumber\\
&+\int_0^t \int_{\mathbb{R}} c(s,x,Y(s,x,z),u(s,x,z),z,\zeta)_{z=Z}\tilde{N}(ds,d\zeta)\nonumber\\
&= \int_{\mathbb{R}} \xi(x,z)\delta_Z(z)dz+\int_0^t \int_{\mathbb{R}}[A_uY(s,x,z)+a(s,x,Y(s,x,z),u(s,x,z),z)]\delta_Z(z)dzds\nonumber\\
&+ \int_0^t \int_{\mathbb{R}}b(s,x,Y(s,x,z),u(s,x,z),z)\delta_Z(z)dzdB(s)\nonumber\\
&+\int_0^t \int_{\mathbb{R}}\int_{\mathbb{R}} c(s,x,Y(s,x,z),u(s,x,z),z,\zeta)\delta_Z(z)dz\tilde{N}(ds,d\zeta)\nonumber\\
&=\int_{\mathbb{R}} \{\xi(x,z)+\int_0^t[A_uY(s,x,z)+ a(s,x,Y(s,x,z),u(s,x,z),z)]ds + \int_0^t b(s,x,Y(s,x,z),u(s,x,z),z)dB(s)\nonumber\\
&+\int_0^t \int_{\mathbb{R}} c(s,x,Y(s,x,z),u(s,x,z),z,\zeta)\tilde{N}(ds,d\zeta)\} \delta_Z(z)dz.
\end{align}
Comparing \eqref{eq1.6} and \eqref{eq1.7} we see that  \eqref{eq1.6} holds if we for each $z$ choose $Y(t,x,z)$ as the solution of the classical (but parametrized) SPDE
\begin{equation} \label{eq3.3}
\begin{cases}
dY(t,x,z) = [A_uY(t,x,z)+a(t,x,Y(t,x,z),u(t,x,z),z)]dt + b(t,x,Y(t,x,z),u(t,x,z),z)dB(t)\\
    + \int_{\mathbb{R}} c(t,x,Y(t,x,z),u(t,x,z),z,\zeta)\tilde{N}(dt,d\zeta); \quad (t,x)\in (0,T)\times D\\
   Y(0,x,z)=\xi(x,z); \quad x\in D \\
\int_{\mathbb{R}}Y(t,x,z)\delta_Z(z)dz=
\theta(t,x);
\quad  (t,x)\in[0,T]\times \partial D.
\end{cases}
\end{equation}

As before let $\mathcal{A}$ be the given family of admissible $\mathbb{H}-$adapted controls $u$.
Then in terms of $Y(t,x,z)$ the performance functional $J(u)$ of a control process $u \in \mathcal{A}$ defined in \eqref{eq1.4} gets the form
\begin{align} \label{eq0.13}
  J(u) &= \mathbb{E}[\int_0^T (\int_Dh(t,x, Y(t,x,Z),u(t,x,Z),Z)dx)dt +\int_Dk(x,Y(T,x,Z),Z)dx] \nonumber\\
   &= \mathbb{E}[\int_\mathbb{R}\Big\{ \int_0^T (\int_{D}h(t, Y(t,x,z),u(t,x,z),z)\mathbb{E}[\delta_Z(z)|\mathcal{F}_t]dx)dt \nonumber \\
   & +\int_D k(x,Y(T,x,z),z)\mathbb{E}[\delta_Z(z)|\mathcal{F}_T]dx\Big\}dz]\nonumber\\
   &= \int_{\mathbb{R}} j(u)(z) dz,
\end{align}
where
\begin{align}\label{eq1.5}
j(u)(z)&:=  \mathbb{E}[ \int_0^T (\int_{D}h(t, Y(t,x,z),u(t,x,z),z)\mathbb{E}[\delta_Z(z)|\mathcal{F}_t]dx)dt \nonumber \\
   & +\int_D k(x,Y(T,x,z),z)\mathbb{E}[\delta_Z(z)|\mathcal{F}_T]dx.
 \end{align}

Thus we see that to maximize $J(u)$ it suffices to maximize $j(u)(z)$ for each value of the parameter $z \in \mathbb{R}$. Therefore Problem 1.1 is transformed into the problem
\begin{problem}
For each given $z \in \mathbb{R}$ find $u^{\star}=u^{\star}(t,x,z) \in\mathcal{A}$ such that
\begin{equation}\label{problem2}
    \sup_{u\in\mathcal{A}}j(u)(z)=j(u^{\star})(z).
\end{equation}
\end{problem}

\section{A sufficient-type maximum principle}
In this section we will establish a sufficient maximum principle for Problem 3.1.\\
We first recall some basic concepts and results from Banach space theory.\\
 Let $\mathcal{X}$ be a Banach space with norm $\| \cdot \|$ and let $F : \mathcal{X} \rightarrow \RB$.
\begin{myenumerate}
\item We say that $F$ has a directional derivative (or G\^{a}teaux derivative) at $v \in \mathcal{X}$ in the direction $w \in \mathcal{X}$ if
$$D_w F(v) := \lim_{\varepsilon \rightarrow 0} \frac{1}{\varepsilon} (F(v + \varepsilon w) - F(v))$$ exists.
\item We say that $F$ is  Fr\'{e}chet differentiable at $v \in V$ if there exists a
 continuous   linear map
$A: \mathcal{X} \rightarrow \RB$
such that
$$\lim_{\substack{h \rightarrow 0 \\ h \in \mathcal{X}}} \frac{1}{ \|h\|} | F(v+h) - F(v) - A(h)| = 0.$$
In this case we call $A$ the {\it gradient} (or Fr\'{e}chet derivative) of $F$ at $v$ and we write
$$A =\nabla_v F.$$
\item If $F$ is Fr\'{e}chet differentiable, then $F$ has a directional derivative in all directions $w \in \mathcal{X}$ and
$$D_w F(v) = \nabla_v F(w) =: \langle \nabla_vF,w\rangle .$$
\end{myenumerate}
In particular, note that if $F$ is a linear operator, then $ \nabla_vF=F$ for all $v$.\\


Problem 3.1 is a stochastic control problem with a standard (albeit parametrized) stochastic partial differential equation \eqref{eq3.3} for the state process $Y(t,x,z)$, but with a non-standard performance functional given by \eqref{eq1.5}. We can solve this problem by a modified maximum principle approach, as follows:\\

Define the \emph{Hamiltonian}
 $ H:[0,T]\times D\times \mathbb{R}\times \mathcal{D}\times\mathbb{U}\times\mathbb{R}\times\mathbb{R}\times\mathbb{R}\times\mathcal{R} \times \Omega \rightarrow \mathbb{R}$ by
\begin{align}\label{eq4.1}
&H(t,x,y,\varphi,u,z,p,q,r)=H(t,x,y,\varphi,u,z,p,q,r,\omega)\nonumber\\
&=\mathbb{E}[\delta_Z(z)|\mathcal{F}_t] h(t,x,y,u,z)+[A_u(\varphi)+a(t,x,y,u,z)]p\nonumber\\
& + b(t,x,y,u,z)q+\int_{\mathbb{R}}c(t,x,y,u,z,\zeta)r(\zeta)\nu(d\zeta).
\end{align}
Here $\mathcal{D}$ denotes the domain of definition for the operator $A_u$, while $\mathcal{R}$ denotes the set of all functions $r(\cdot) : \mathbb{R}\rightarrow  \mathbb{R}$
such that the last integral above converges. We assume that $\mathcal{D}$ is a Banach space.The quantities $p,q,r(\cdot)$ are called the \emph{adjoint variables}.
The \emph{adjoint processes} $p(t,x,z),q(t,x,z),r(t,x,z,\zeta)$ are defined as the solution of the $z$-parametrized backward stochastic partial differential equation (BSPDE)
\begin{equation}\label{eq4.2a}
    \left\{
\begin{array}{l}
    dp(t,x,z) = -[A_{u(t,x,z)}^{*}p(t,x,z) + \frac{\partial H}{\partial y}(t,x,z)]dt +q(t,x,z)dB(t)+\int_{\mathbb{R}}r(t,x,z,\zeta)\tilde{N}(dt,d\zeta); \\ \quad (t,x,z)\in(0,T)\times D\times \mathbb{R}\\
    p(T,x,z)  = \frac{\partial k}{\partial y}(x,Y(T,x,z),z) \mathbb{E}[\delta_Z(z)|\mathcal{F}_T];\quad (x,z)\in D\times \mathbb{R}\\
    p(t,x,z)=0; \quad (t,x,z)\in [0,T]\times \partial D\times\mathbb{R},
\end{array}
    \right.
\end{equation}
where
\begin{equation}
\frac{\partial H}{\partial y}(t,x,z)=\frac{\partial H}{\partial y}(t,x,y,Y(t,.,z),u(t,x,z),z,p(t,x,z),q(t,x,z),r(t,x,z,.))|_{y=Y(t,x,z)}.
\end{equation}
For fixed $t,u,z,p,q,r$ we can regard
\begin{equation}\label{eq4.5a}
\varphi \mapsto \ell(\varphi)(x):=H(t,x,\varphi(x),\varphi,u,z,p,q,r)
\end{equation}
as a map from $\mathcal{D}$ into $\mathbb{R}$. The Fr\'{e}chet derivative at $\varphi$ of this map is the linear operator $\nabla_{\varphi}\ell$ on $\mathcal{D}$ given by
\begin{equation}
\langle \nabla_{\varphi}\ell,\psi \rangle=\nabla_{\varphi}\ell(\psi)=A_u(\psi)(x)p+ \psi(x)\frac{\partial H}{\partial y}(t,x,y,\varphi,u,z,p,q,r)|_{y=\varphi(x)}; \quad \psi \in \mathcal{D}.
\end{equation}
For simplicity of notation, if there is no risk of confusion,  we will denote $\ell$ by $H$ from now on.


We can now state the first maximum principle for our problem \eqref{problem2}:

\begin{theorem}{[Sufficient-type maximum principle]}\\
Let $\hat{u} \in \mathcal{A}$, and denote the associated solution  of \eqref{eq3.3} and \eqref{eq4.2a} by $\hat{Y}(t,x,z)$ and \\
$(\hat{p}(t,x,z),\hat{q}(t,x,z),\hat{r}(t,x,z,\zeta))$, respectively.  Assume that the following hold:
\begin{enumerate}
 \item $ y \rightarrow k(x,y,z)$ is concave for all $x,z$
 \item $(\varphi,u)\rightarrow H(t,x,\varphi(x), \varphi,u,z,\widehat{p}(t,x,z),\widehat{q}(t,x,z),\hat{r}(t,x,z,\zeta))$ is concave for all $t,x,z,\zeta$
 \item $\sup_{w\in\mathbb{U}}H\big(t,x,\widehat{Y}(t,x,z),\widehat{Y}(t,\cdot,z)(x),w,\widehat{p}(t,x,z),\widehat{q}(t,x,z),\hat{r}(t,x,z,\zeta)\big)$\\
      $=H\big(t,x,\widehat{Y}(t,x,z), \widehat{Y}(t,\cdot,z)(x),\widehat{u}(t,x,z),\widehat{p}(t,x,z),\widehat{q}(t,x,z),\hat{r}(t,x,z,\zeta)\big)$ for all $t,x,z,\zeta.$
 \end{enumerate}
Then $\widehat{u}(\cdot,\cdot,z)$ is an optimal insider control for Problem 3.1.
\end{theorem}
\dproof  By considering an increasing sequence of stopping times $\tau_n$ converging to $T$, we may assume that all local integrals appearing in the computations below are martingales and hence have expectation 0. See \cite{OS2}. We omit the details.\\
Choose arbitrary $u(.,.,z)\in\mathcal{A}$, and let the corresponding solution of \eqref{eq3.3} and \eqref{eq4.2a}  be $Y(t,x,z)$, $p(t,x,z)$, $q(t,x,z)$, $r(t,x,z,\zeta)$.
For simplicity of notation we write\\
$h=h(t,x,Y(t,x,z),u(t,x,z))$,
$\widehat{h}=h(t,x,\widehat{Y}(t,x,z),\widehat{u}(t,x,z))$
 and similarly with $a$, $\widehat{a}$, $b$, $\widehat{b}$ and so on.\\
 Moreover put
 \begin{equation}
 \hat{H}(t,x)=H(t,x,\widehat{Y}(t,x,z),\widehat{Y}(t,\cdot,z)(x),\widehat{u}(t,x,z),\widehat{p}(t,x,z),\widehat{q}(t,x,z),\widehat{r}(t,x,z,.))
 \end{equation}
 and
 \begin{equation}
 H(t,x)=H(t,x,Y(t,x,z), Y(t,\cdot,z)(x),u(t,x,z),\widehat{p}(t,x,z),\widehat{q}(t,x,z),\widehat{r}(t,x,z,.))
 \end{equation}
  In the following we write  $\widetilde{h}=h-\widehat{h}$, $\widetilde{a}=a-\widehat{a}$, $\widetilde{Y}=Y-\widehat{Y}$.\\
 Consider
 \begin{equation*}
    j(u(.,.,z))-j(\widehat{u}(.,.,z))=I_1+ I_2,
 \end{equation*}
 where
 \begin{equation}\label{eq4.7}
    I_1=\mathbb{E}[\int_0^T(\int_D\{h(t,x)-\widehat{h}(t,x)\}\mathbb{E}[\delta_Z(z)|\mathcal{F}_t]dx)dt], \quad I_2=\mathbb{E}[\int_D\{k(x)-\hat{k}(x)\}\mathbb{E}[\delta_Z(z)|\mathcal{F}_T]dx].
 \end{equation}
  By the definition of $H$ we have
  \begin{align}\label{eq4.8}
    I_1 &= \mathbb{E}[\int_0^T\int_D\{H(t,x)-\widehat{H}(t,x)-\widehat{p}(t,x) [A_uY(t,x)-A_{\hat{u}}\widehat{Y}(t,x) +\widetilde{a}(t,x)] - \widehat{q}(t,x)\widetilde{b}(t,x)\nonumber\\
   & -\int_{\mathbb{R}}\hat{r}(t,x,\zeta)\tilde{c}(t,x,\zeta)\nu(d\zeta)\}dx dt].
  \end{align}

 Since $k$ is concave with respect to $y$ we have
\begin{align}
&(k(x,Y(T,x,z),z)-k(x,\hat{Y}(T,x,z),z))\mathbb{E}[\delta_Z(z)|\mathcal{F}_T]\nonumber\\
&\leq \frac{\partial k}{\partial y}(x,\hat{Y}(T,x,z),z)\mathbb{E}[\delta_Z(z)|\mathcal{F}_T](Y(T,x,z)-\hat{Y}(T,x,z)),
\end{align}
and hence
 \begin{align}\label{eq4.11}
   I_2 &\leq\mathbb{E}[\int_D\frac{\partial k}{\partial y}(x,\widehat{Y}(T,x,z))\mathbb{E}[\delta_Z(z)|\mathcal{F}_T]\tilde{Y}(T,x,z)dx]=\mathbb{E}[\int_D\widehat{p}(T,y)\widetilde{Y}(T,x,z)dx] \\ \nonumber
    &= \mathbb{E}[\int_D(\int_0^T \widehat{p}(t,x,z) d\widetilde{Y}(t,x,z)+\int_0^T\widetilde{Y}(t,x,z)d\widehat{p}(t,x,z)+\int_0^Td[\hat{p} , \tilde{Y}]_t)dx] \\ \nonumber
    &= \mathbb{E}[\int_D\int_0^T \{\widehat{p}(t,x,z) [A_uY(t,x,z)-A_{\hat{u}}\widehat{Y}(t,x)+\widetilde{a}(t,x)] -\widetilde{Y}(t,x,z) [A_{\hat{u}}^{*}\widehat{p}(t,x,z)+\frac{\partial \widehat{H}(t,x)}{\partial y}]\nonumber\\
   & + \widetilde{b}(t,x)\widehat{q}(t,x)+\int_{\mathbb{R}}\tilde{c}(t,x,z,\zeta)\hat{r}(t,x,z,\zeta)\nu(d\zeta)\}dtdx].\nonumber
    \end{align}
where
\begin{equation}
\frac{\partial \widehat{H}(t,x)}{\partial y}=\frac{\partial H}{\partial y}(t,x,\hat{Y}(t,x,z),\widehat{Y}(t,\cdot,z)(x),\hat{u}(t,x,z),\hat{p}(t,x,z),\hat{q}(t,x,z),\hat{r}(t,x,z,.)).
\end{equation}
By a slight extension of \eqref{adjoint} we get
\begin{equation}\label{eq4.15a}
\int_{D} \widetilde{Y}(t,x,z)A_{\hat{u}}^{*}\widehat{p}(t,x,z)dx =\int_{D} \widehat{p}(t,x,z)A_{\hat{u}}\widetilde{Y}(t,x,z)dx.
\end{equation}

Therefore, adding \eqref{eq4.8} - \eqref{eq4.11} and using \eqref{eq4.15a} we get,
\begin{align}\label{eq4.10aa}
j(u(.,z))-j(\widehat{u}(.,z))
\leq\mathbb{E}[\int_{D}\Big(\int_0^T\{H(t,x)-\hat{H}(t,x)-[\hat{p}(t,x,z)A_{\hat{u}}(\tilde{Y})(t,x,z)+\tilde{Y}(t,x,z)\frac{\partial \hat{H}(t,x)}{\partial y}]\}dt \Big)dx].
\end{align}

Hence
\begin{align}\label{eq4.10}
j(u(.,z))-j(\widehat{u}(.,z))
\leq\mathbb{E}[\int_{D}\Big(\int_0^T\{H(t,x)-\hat{H}(t,x)-\nabla_{\hat{Y}} \widehat{H}(\tilde{Y})(t,x,z)\}dt \Big)dx]
\end{align}
where
\begin{equation}
\nabla_{\hat{Y}} \widehat{H}(\tilde{Y})=\nabla_{\varphi} \widehat{H}(\tilde{Y})|_{\varphi=\hat{Y}}
\end{equation}
By the concavity assumption of $H$ in $(\varphi,u)$ we have:
\begin{equation}
H(t,x)-\hat{H}(t,x)\leq \nabla_{\hat{Y}} \widehat{H}(Y-\hat{Y})(t,x,z)+ \frac{\partial \widehat{H}}{\partial u}(t,x)(u(t,x)-\hat{u}(t,x)),
\end{equation}
and the maximum condition implies that
\begin{equation}
\frac{\partial\widehat{H}}{\partial u}(t,x)(u(t,x)-\hat{u}(t,x))\leq 0.
\end{equation}
Hence by \eqref{eq4.10} we get
$j(u)\leq j(\hat{u})$.
Since $u\in\mathcal{A}$ was arbitrary, this shows that $\hat{u}$ is optimal.

\fproof

\section{A necessary-type maximum principle}

We proceed to establish a corresponding necessary maximum principle. For this, we do not need concavity conditions, but instead we need the following assumptions about the set of admissible control processes:\\
\begin{itemize}
\item
$A_1$. For all $t_0\in [0,T]$ and all bounded $\mathcal{H}_{t_0}$-measurable random variables $\alpha(x,z,\omega)$, the control
$\theta(t,x,z, \omega) := \mathbf{1}_{[t_0,T ]}(t)\alpha(x,z,\omega)$ belongs to $\mathcal{A}$.\\
\item
$A_2$. For all $u, \beta_0 \in\mathcal{A}$ with $\beta_0(t,x,z) \leq K < \infty$ for all $t,x,z$  define
\begin{equation}\label{delta}
    \delta(t,x,z)=\frac{1}{2K}dist(u(t,x,z),\partial\mathbb{U})\wedge1 > 0
\end{equation}
and put
\begin{equation}\label{eq3.2}
    \beta(t,x,z)=\delta(t,x,z)\beta_0(t,x,z).
\end{equation}
Then the control
\begin{equation*}
    \widetilde{u}(t,x,z)=u(t,x,z) + a\beta(t,x,z) ; \quad t \in [0,T]
\end{equation*}
belongs to $\mathcal{A}$ for all $a \in (-1, 1)$.\\
\item
$A3$. For all $\beta$ as in \eqref{eq3.2} the derivative process
\begin{equation}\label{eq5.3a}
    \chi(t,x,z):=\frac{d}{da}Y^{u+a\beta}(t,x,z)|_{a=0}
\end{equation}
exists, and belong to $\mathbf{L}^2(\lambda\times \mathbf{P})$ and

\begin{equation}\label{d chi}
    \left\{
\begin{array}{l}
    d\chi(t,x,z) = [\frac{dA}{du}(Y)(t,x,z)\beta(t,x,z)+A_u \chi(t,x,z)+\frac{\partial a}{\partial y}(t,x,z)\chi(t,x,z)\\
    +\frac{\partial a}{\partial u}(t,x,z)\beta(t,x,z)]dt\\
    +[\frac{\partial b}{\partial y}(t,x,z)\chi(t,x,z)+\frac{\partial b}{\partial u}(t,x,z)\beta(t,x,z)]d B(t) \\
    +\int_{\mathbb{R}}[\frac{\partial c}{\partial y}(t,x,z,\zeta)\chi(t,x,z)+\frac{\partial c}{\partial u}(t,x,z,\zeta)\beta(t,x,z)]\tilde{N}(dt,d\zeta); (t,x) \in [0,T]\times D,\\
    \chi(0,x,z)  = \frac{d}{da}Y^{u+a\beta}(0,x,z)|_{a=0} = 0,\\
    \chi(t,x,z) = 0; (t,x) \in [0,T]\times \partial D.
   \end{array}
    \right.
\end{equation}
\end{itemize}

\begin{theorem}{[Necessary-type maximum principle]}\label{necessary theorem} \\
Let $\hat{u} \in \mathcal{A}$ and $z \in \mathbb{R}$. Then the following are equivalent:
\begin{enumerate}
\item $\frac{d}{da}j(\hat{u}+a\beta)(z)|_{a=0}=0$ for all bounded $\beta \in \mathcal{A}$ of the form \eqref{eq3.2}.
\item $\frac{\partial H}{\partial u}(t,x,z)_{u=\hat{u}}=0$ for all $(t,x)\in [0,T] \times D.$
\end{enumerate}
\end{theorem}

\dproof
For simplicity of notation we write $u$ instead of $\hat{u}$ in the following. \\
By considering an increasing sequence of stopping times $\tau_n$ converging to $T$, we may assume that all local integrals appearing in the computations below are martingales and have expectation 0. See \cite{OS2}. We omit the details.\\
We can write
$$\frac{d}{da}j((u+a\beta)(z))|_{a=0}=I_1+I_2$$\\
where
$$I_1=\frac{d}{da}\mathbb{E}[\int_D\int_0^Th(t,x,Y^{u+a\beta}(t,x,z),u(t,x,z)+a\beta(t,x,z),z)\mathbb{E}[\delta_Z(z)|\mathcal{F}_t]dtdx]|_{a=0}$$\\
and
$$I_2=\frac{d}{da}\mathbb{E}[\int_Dk(x,Y^{u+a\beta}(T,x,z),z)\mathbb{E}[\delta_Z(z)|\mathcal{F}_T]dx]|_{a=0}.$$
By our assumptions on $h$ and $k$ and by \eqref{eq5.3a} we have
\begin{equation}\label{iii1}
    I_1=\mathbb{E}[\int_D\int_0^T\{\frac{\partial h}{\partial y}(t,x,z)\chi(t,x,z)+\frac{\partial h}{\partial u}(t,x,z)\beta(t,x,z)\}\mathbb{E}[\delta_Z(z)|\mathcal{F}_t]dtdx],
\end{equation}
\begin{equation}\label{iii2}
    I_2=\mathbb{E}[\int_D\frac{\partial k}{\partial y}(x,Y(T,x,z),z)\chi(T,x,z)\mathbb{E}[\delta_Z(z)|\mathcal{F}_T]dx]=\mathbb{E}[\int_Dp(T,x,z)\chi(T,x,z)dx].
\end{equation}
By the It\^{o} formula

\begin{align}\label{eq5.9a}
   I_2&= \mathbb{E}[\int_Dp(T,x,z)\chi(T,x,z)dx]=\mathbb{E}[\int_D\int_0^Tp(t,x,z)d\chi(t,x,z) dx+\int_D\int_0^T\chi(t,x,z)dp(t,x,z) dx\nonumber\\
   &+\int_D\int_0^Td[\chi,p](t,x,z)dx] \\ \nonumber
   &= \mathbb{E}[\int_D\int_0^Tp(t,x,z)\{\frac{dA}{du}(Y)(t,x,z)\beta(t,x,z)+A_u \chi(t,x,z)\nonumber\\
   &+\frac{\partial a}{\partial y}(t,x,z)\chi(t,x,z)+\frac{\partial a}{\partial u}(t,x,z)\beta(t,x,z)\}dtdx\\ \nonumber
   &+\int_D\int_0^Tp(t,x,z)\{\frac{\partial b}{\partial y}(t,x,z)\chi(t,x,z)+\frac{\partial b}{\partial u}(t,x,z)\beta(t,x,z)\}dB(t) \\ \nonumber
   &+\int_D\int_0^T\int_{\mathbb{R}}p(t,x,z)\{\frac{\partial c}{\partial y}(t,x,z,\zeta)\chi(t,x,z)+\frac{\partial c}{\partial u}(t,x,z,\zeta)\beta(t,x,z)\}\tilde{N}(dt,d\zeta)dx\\ \nonumber
   &-\int_D\int_0^T\chi(t,x,z)[A_u^*p(t,x,z)+\frac{\partial H}{\partial y}(t,x,z)] dt dx \nonumber\\
   &+\int_D\int_0^T\chi(t,x,z)q(t,x,z)dB(t)dx+\int_D\int_0^T\int_{\mathbb{R}}\chi(t,x,z)r(t,x,z,\zeta) \tilde{N}(dt,d\zeta)dx  \\ \nonumber
     &+\int_D\int_0^Tq(t,x,z) \{\frac{\partial b}{\partial y}(t,x,z)\chi(t,x,z)+\frac{\partial b}{\partial u}(t,x,z)\beta(t,x,z)\}dt dx \\ \nonumber
   &+\int_D\int_0^T\int_{\mathbb{R}}\{\frac{\partial c}{\partial y}(t,x,z,\zeta)\chi(t,x,z)+\frac{\partial c}{\partial u}(t,x,z,\zeta)\beta(t,x,z)\}r(t,x,z,\zeta)\nu(\zeta)dt dx]\\ \nonumber
  &= \mathbb{E}\Big[\int_D\Big(\int_0^T\{p(t,x,z)(\frac{dA}{du}(Y)(t,x,z)\beta(t,x,z)+A_u\chi(t,x,z))\}dt\nonumber\\
  &+\int_0^T\chi(t,x,z)\{p(t,x,z)\frac{\partial a}{\partial y}(t,x,z)+q(t,x,z)\frac{\partial b}{\partial y}(t,x,z)-A_u^*p(t,x,z)-\frac{\partial H}{\partial y}(t,x,z)\nonumber\\
  &+\int_{\mathbb{R}}\frac{\partial c}{\partial y}(t,x,z,\zeta)r(t,x,z,\zeta)\nu(d\zeta)\}dt  \\ \nonumber
   &+ \int_0^T\beta(t,x,z)\{p(t,x,z)\frac{\partial a}{\partial u}(t,x,z)+q(t,x,z)\frac{\partial b}{\partial u}(t,x,z)+\int_{\mathbb{R}}\frac{\partial c}{\partial u}(t,x,z,\zeta)r(t,x,z,\zeta)\nu(d\zeta)\}dt\Big)dx\Big]\\ \nonumber
&=\mathbb{E}\Big[\int_D[\int_0^T -\chi(t,x,z)\frac{\partial h}{\partial y}\mathbb{E}[\delta_Z(z)|\mathcal{F}_t]\}dt\nonumber\\
&+\int_0^T\{\frac{\partial H}{\partial u}(t,x,z)-\frac{\partial h}{\partial u}(t,x,z)\mathbb{E}[\delta_Z(z)|\mathcal{F}_t]\}\beta(t,x,z)dt]dx\Big] \\
&= -I_1+\mathbb{E}[\int_D\int_0^T\frac{\partial H}{\partial u}(t,x,z)\beta(t,x,z)dt dx].\nonumber
\end{align}

Summing (\ref{iii1}) and \eqref{eq5.9a} we get
\begin{equation*}
    \frac{d}{da}j((u+a\beta)(.,x,y))|_{a=0}=I_1+I_2=\mathbb{E}[\int_D\int_0^T\frac{\partial H}{\partial u}(t,x,z)\beta(t,x,z)dtdx].
\end{equation*}
We conclude that
\begin{equation*}
    \frac{d}{da}j(u+a\beta)(z))|_{a=0}=0
\end{equation*}
if and only if
\begin{equation}
\mathbb{E}[\int_D\int_0^T\frac{\partial H}{\partial u}(t,x,z)\beta(t,x,z)dtdx]=0,
\end{equation}
for all bounded $\beta\in\mathcal{A}$ of the form \eqref{eq3.2}.\\

\noindent In particular, applying this to $\beta(t,x,z) = \theta(t,x,z)$ as in $A1$, we get that this is again equivalent to
\begin{equation}
   \frac{\partial H}{\partial u}(t,x,z)=0  \text{ for all } (t,x) \in [0,T] \times D.
\end{equation}
\fproof
\section{Controls which do not depend on $x$}
In some situations it is of interest to study controls $u(t, x)=u(t)$ which have the same value throughout the space $\mathcal{D}$, i.e., only depends on time $t$. See e.g. Section 8.2.  In this case we define the set $\mathcal{A}_0$ of admissible
controls by
\begin{equation}
\mathcal{A}_0 = \{u\in\mathcal{A}; u(t, x) = u(t) \text{ does not depend on } x\}.
\end{equation}
Defining the performance functional $J(u)=\int_{\mathbb{R}}j(u)(z)dz$ as in Problem 3.1, the problem now becomes:
\begin{problem}
For each $z \in \mathbb{R}$ find  $u_0^{\ast}\in\mathcal{A}_0$ such that
\begin{equation}\label{problem4}
sup_{u\in\mathcal{A}_0}j(u)(z) = j(u_0^{\ast})(z).
\end{equation}
\end{problem}

\subsection{Sufficient-type maximum principle for controls which do not depend on $x$}
We now state and prove an analog of Theorem 4.1 for this case:
\begin{theorem}(Sufficient-type maximum principle for controls which do not depend on x).
Suppose $\hat{u}\in\mathcal{A}_0$ with corresponding solutions $\hat{Y}(t, x, z)$ of \eqref{eq3.3}
and $\hat{p}(t, x, z), \hat{q}(t, x, z), \hat{r}(t, x, z, \zeta)$ of \eqref{eq4.2a}
respectively. Assume that the following hold:\\
\begin{enumerate}
 \item $ y \rightarrow k(x,y,z)$ is concave for all $x,z$
 \item $(\varphi,u)\rightarrow H(t,x,\varphi(x),\varphi,u,z,\widehat{p}(t,x,z),\widehat{q}(t,x,z),\hat{r}(t,x,z,\cdot))$ is concave for all $t,x,z$
 \item $\sup_{w\in\mathbb{U}}\int_D H\big(t,x,\widehat{Y}(t,x,z),\widehat{Y}(t,\cdot,z), w,\widehat{p}(t,x,z),\widehat{q}(t,x,z),\hat{r}(t,x,z,\cdot)\big)dx$\\
      $=\int_D H\big(t,x,\widehat{Y}(t,x,z),\widehat{Y}(t,\cdot,z),\widehat{u}(t,z),\widehat{p}(t,x,z),\widehat{q}(t,x,z),\hat{r}(t,x,z,\cdot)\big)dx$ for all $t,z.$
 \end{enumerate}
Then $\hat{u}(t,z)$ is an optimal control for the Problem 6.1.
\end{theorem}

\dproof
We proceed as in the proof of Theorem 4.1. Let $u\in \mathcal{A}_0$ with corresponding
solution $Y(t,x,z)$ of \eqref{eq3.3}. With $\hat{u}\in\mathcal{A}_0$,
consider
\begin{equation}\label{eq4.7}
j(u) - j(\hat{u}) = \mathbb{E}[\int_0^T\int_D \{h - \hat{h}\}dx dt +\int_D \{k - \hat{k}\}dx] ,
\end{equation}
where\\

$\hat{h} = h(t, x,\hat{Y} (t,x,z), \hat{u}(t,z)),\quad h = h(t, x, Y(t, x,z), u(t,z))$\\

$\hat{k} = k(x, \hat{Y} (T, x,z)) \text{  and  } k = k(x, Y(T, x,z))$.

Using a similar shorthand notation for  $\hat{a}, a, \hat{b}, b$ and  $\hat{c}, c$,  and setting
\begin{equation}
\hat{H}= H(t, x,\hat{Y} (t, x,z), \hat{u}(t,z), \hat{p}(t, x,z), \hat{q}(t, x,z), \hat{r}(t,x,z,\cdot))
\end{equation}
and
\begin{equation}
H = H(t, x, Y(t, x,z), u(t,z), \hat{p}(t, x), \hat{q}(t, x), \hat{r}(t, x,z, \cdot)),
\end{equation}
we see that (\ref{eq4.7}) can be written
\begin{equation}
j(u) - j(\hat{u}) = I_1 + I_2,
\end{equation}

where
 \begin{equation}\label{I_1I_2}
    I_1=\mathbb{E}[\int_0^T(\int_D\{h(t,x)-\widehat{h}(t,x)\}\mathbb{E}[\delta_Z(z)|\mathcal{F}_t]dx)dt], \quad I_2=\mathbb{E}[\int_D\{k(x)-\hat{k}(x)\}\mathbb{E}[\delta_Z(z)|\mathcal{F}_T]dx].
 \end{equation}
  By the definition of $H$ we have
  \begin{align}\label{II1}
    I_1 &= \mathbb{E}[\int_0^T\int_D\{H(t,x)-\widehat{H}(t,x)-\widehat{p}(t,x)(A_uY(t,x)-A_{\hat{u}}\hat{Y}(t,x)+\widetilde{a}(t,x)) - \widehat{q}(t,x)\widetilde{b}(t,x)\nonumber\\
   & -\int_{\mathbb{R}}\hat{r}(t,x,\zeta)\tilde{c}(t,x,\zeta)\nu(d\zeta)\}dx dt].
  \end{align}

 Since $k$ is concave with respect to $y$ we have
\begin{align}
&(k(x,Y(T,x,z),z)-k(x,\hat{Y}(T,x,z),z))\mathbb{E}[\delta_Z(z)|\mathcal{F}_T]\nonumber\\
&\leq \frac{\partial k}{\partial y}(x,\hat{Y}(T,x,z),z)\mathbb{E}[\delta_Z(z)|\mathcal{F}_T](Y(T,x,z)-\hat{Y}(T,x,z)).
\end{align}
Therefore, as in the proof of Theorem 4.1,
 \begin{align}\label{II_2}
   I_2 &\leq\mathbb{E}[\int_D\int_0^T \{\widehat{p}(t,x,z) [A_uY(t,x,z)-A_{\tilde{u}}\tilde{Y}(t,x,z)+\widetilde{a}(t,x)] -\widetilde{Y}(t,x,z) [A_{\hat{u}}^*\hat{p}(t,x,z)
   +\frac{\partial \widehat{H}(t,x)}{\partial y}(t,x)]\nonumber\\
   & + \widetilde{b}(t,x)\widehat{q}(t,x)+\int_{\mathbb{R}}\tilde{c}(t,x,z,\zeta)\hat{r}(t,x,z,\zeta)\nu(d\zeta)\}dtdx].
    \end{align}
where
\begin{equation}
\frac{\partial \widehat{H}(t,x)}{\partial y}=\frac{\partial H}{\partial y}(t,x,\hat{Y}(t,x,z),\hat{Y}(t,.,z)(x),\hat{u}(t,z),\hat{p}(t,x,z),\hat{q}(t,x,z),\hat{r}(t,x,z,.))
\end{equation}

Adding (\ref{II1}) - (\ref{II_2}) we get as in equation \eqref{eq4.10},
\begin{equation} \label{eq4.17}
j(u)-j(\hat{u})\leq \mathbb{E}\Big[\int_0^T \Big( \int_D \{H(t,x)-\hat{H}(t,x)-\nabla_{\hat{Y}}\hat{H}(\tilde{Y})(t,x,z)\}dx \Big)dt \Big].
\end{equation}
By the concavity assumption of $H$ in $(y,u)$ we have:
\begin{equation}
H(t,x)-\hat{H}(t,x)\leq \nabla_{\hat{Y}}\hat{H}(Y-\hat{Y})(t,x,z) + \frac{\partial \widehat{H}}{\partial u}(t,x)(u(t)-\hat{u}(t)),
\end{equation}
and the maximum condition implies that
\begin{equation}
\int_D \frac{\partial\widehat{H}}{\partial u}(t,x)(u(t)-\hat{u}(t)) dx\leq 0.
\end{equation}
Hence
\begin{equation}
\int_D \{ H(t,x)-\hat{H}(t,x) -\nabla_{\hat{Y}}\hat{H}(Y-\hat{Y})(t,x,z) \}dx \leq 0,
\end{equation}
and therefore we conclude by \eqref{eq4.17} that $j(u)\leq j(\hat{u})$.
Since $u\in\mathcal{A}_0$ was arbitrary, this shows that $\hat{u}$ is optimal.
\fproof
\subsection{Necessary-type maximum principle for controls which do not depend on $x$}
We proceed as in Theorem \ref{necessary theorem} to establish a corresponding necessary maximum principle for controls which do not depend on $x$. As in Section 5 we assume the following:

\begin{itemize}
\item
$A_1$. For all $t_0\in [0,T]$ and all bounded $\mathcal{H}_{t_0}$-measurable random variables $\alpha(z,\omega)$, the control
$\theta(t,z, \omega) := \mathbf{1}_{[t_0,T ]}(t)\alpha(z,\omega)$ belongs to $\mathcal{A}_0$.\\
\item
$A_2$. For all $u, \beta_0 \in\mathcal{A}_0$ with $\beta_0(t,z) \leq K < \infty$ for all $t,z$  define
\begin{equation}\label{delta}
    \delta(t,z)=\frac{1}{2K}dist((u(t,z),\partial\mathbb{U})\wedge1 > 0
\end{equation}
and put
\begin{equation}\label{eq4.22}
    \beta(t,z)=\delta(t,z)\beta_0(t,z).
\end{equation}
Then the control
\begin{equation*}
    \widetilde{u}(t,z)=u(t,z) + a\beta(t,z) ; \quad t \in [0,T]
\end{equation*}
belongs to $\mathcal{A}_0$ for all $a \in (-1, 1)$.\\
\item
$A3$. For all $\beta$ as in (\ref{eq4.22}) the derivative process
\begin{equation*}
    \chi(t,x,z):=\frac{d}{da}Y^{u+a\beta}(t,x,z)|_{a=0}
\end{equation*}
exists, and belong to $\mathbf{L}^2(\lambda\times \mathbf{P})$ and

\begin{equation}\label{d chi}
    \left\{
\begin{array}{l}
    d\chi(t,x,z) = [\frac{dL}{du}(Y)(t,x,z)\beta(t,z)+A_u\chi(t,x,z)\\
    +\frac{\partial a}{\partial y}(t,x,z)\chi(t,x,z)+\frac{\partial a}{\partial u}(t,x,z)\beta(t,z)]dt\\
    +[\frac{\partial b}{\partial y}(t,x,z)\chi(t,x,z)+\frac{\partial b}{\partial u}(t,x,z)\beta(t,z)]d B(t) \\
    +\int_{\mathbb{R}}[\frac{\partial c}{\partial y}(t,x,z,\zeta)\chi(t,x,z)+\frac{\partial c}{\partial u}(t,x,z,\zeta)\beta(t,z)]\tilde{N}(dt,d\zeta); (t,x) \in [0,T] \times D,\\
    \chi(0,x,z)  = \frac{d}{da}Y^{u+a\beta}(0,x,z)|_{a=0} = 0; x \in D\\
    \chi(t,x,z) = 0; (t,x) \in [0,T] \times \partial D.
   \end{array}
    \right.
\end{equation}
\end{itemize}
Then we have the following result:

\begin{theorem}{[Necessary-type maximum principle for controls which do not depend on $x$]}\label{necessary theorem 2} \\
Let $\hat{u} \in \mathcal{A}_0$ and $z \in \mathbb{R}$. Then the following are equivalent:
\begin{enumerate}
\item $\frac{d}{da}j(\hat{u}+a\beta)(z)|_{a=0}=0$ for all bounded $\beta \in \mathcal{A}_0$ of the form \eqref{eq4.22}.
\item $[\int_D\frac{\partial H}{\partial u}(t,x,z)dx]_{u=\hat{u}(t)}=0$ for all $t\in [0,T] .$
\end{enumerate}
\end{theorem}

\dproof
The proof is analogous to the proof of Theorem \ref{necessary theorem} and is omitted.
\fproof
\section{Application to noisy observation optimal control}
For simplicity we consider only the one-dimensional case in the following.\\
Suppose the signal process $X(t) = X^{(u)}(t,Z)$ and its corresponding
observation process $R(t)$ are given respectively by the following system of stochastic differential
equations
\begin{itemize}
\item(Signal process)
\begin{align}\label{eq7.1}
&dX(t) = \alpha(X(t), R(t), u(t,Z))dt + \beta(X(t), R(t), u(t,Z))dv(t)\nonumber\\
&+\int_{\mathbb{R}}\gamma(X(t), R(t), u(t,Z), \zeta) \tilde{N} (dt, d\zeta); t\in[0, T ],\\
&X(0) \text{ has density } F(\cdot), \text{ i.e. } \mathbb{E}[\phi(X(0))] =\int_{\mathbb{R}}\phi(x)F(x)dx;\quad  \phi\in C_0(\mathbb{R}).\nonumber
\end{align}

As before $T > 0$ is a fixed constant.
\item (Observation process)
\begin{align}\label{eq7.2}
\begin{cases}
dR(t) = h(X(t))dt + dw(t); \quad t\in[0, T ],\\
R(0) = 0.
\end{cases}
\end{align}
\end{itemize}
Here $\alpha : \mathbb{R}\times\mathbb{R}\times \mathbb{U} \rightarrow \mathbb{R}, \beta:\mathbb{R}\times\mathbb{R}\times \mathbb{U} \rightarrow \mathbb{R} , \gamma : \mathbb{R}\times\mathbb{R}\times \mathbb{U} \times\mathbb{R} \rightarrow \mathbb{R}$ are given deterministic functions and $h : \mathbb{R}\rightarrow \mathbb{R}$ is a given deterministic function such that the Novikov condition holds, i.e.
\begin{equation}\label{Novikov}
\mathbb{E}[\exp(\frac{1}{2}\int_0^T h^2(X(s,Z))ds)]<\infty.
\end{equation}
 The processes  $v(t) = v(t,\omega)$ and $w(t) = w(t, \omega)$ are independent Brownian motions, and $\tilde{N} (dt, d\zeta)$ is a compensated Poisson random measure, independent of both $v$ and $w$. We let $\mathbb{F}^{v}:= \{ \mathcal{F}^{v}_t\}_{0\leq t \leq T}$ and $ \mathbb{F}^{w}:=\{ \mathcal{F}^{w}_t\}_{0\leq t \leq T}$ denote the filtrations generated by $(v,\tilde{N})$ and $w$, respectively. We assume that $Z$ is a given $\mathcal{F}^{v}_{T_0}$ - measurable random variable, representing the inside information of the controller, where $T_0 > 0$ is a constant. Note that $Z$ is independent of $\mathbb{F}^{w}$.\\

 The process $u(t) =
u(t, Z, \omega)$ is our control process, assumed to have values in a given closed set
$\mathbb{U}\subseteq\mathbb{R}$ . We require that $u(t)$ be adapted to the filtration
\begin{equation}
\mathbb{H}:= \{\mathcal{H}_t\}_{0\leq t \leq T},
\text{ where }
\mathcal{H}_t=\mathcal{R}_t \vee \sigma(Z),
\end{equation}
where $\mathcal{R}_t$ is the sigma-algebra generated by the observations $R(s), s\leq t$.
This means that for all $t$ our control process $u$ is of the form
$$ u=u(t,Z),$$
where $u(t,z)$ is $\mathcal{R}_t$-measurable for each constant $z \in \mathbb{R}$. Similarly the signal process can be written $X=X(t,Z)$, where $X(t,z)$ is the solution of \eqref{eq7.1} with the random variable $Z$ replaced by the parameter $z \in \mathbb{R}$.
We call $u(t)$ admissible if, in addition, \eqref{eq7.1} and \eqref{eq7.2} has a unique strong solution $(X(t), R(t))$ such that
\begin{equation}
\mathbb{E}[\int_0^T|f(X(t), u(t))|dt + |g(X(T ))|]<\infty,
\end{equation}
where $f:\mathbb{R}\times \mathbb{U} \rightarrow\mathbb{R}$ and $g : \mathbb{R}\rightarrow\mathbb{R}$ are given functions, called the profit
rate and the bequest function, respectively. The set of all admissible controls is denoted by $\mathcal{A}_{\mathbb{H}}$.
For $u \in \mathcal{A}_{\mathbb{H}}$ we define the performance functional
\begin{equation}
J(u) = \mathbb{E}[\int_0^Tf(X(t), u(t))dt + g(X(T ))].
\end{equation}
We consider the following problem:\\
\begin{problem}[The noisy observation insider stochastic control problem]
Find
$u^{\ast}\in \mathcal{A}_{\mathbb{H}}$ such that
\begin{equation} \label{eq7.6}
 \sup_{u \in \mathcal{A}_{\mathbb{H}}}J(u) = J(u^{\ast}).
\end{equation}
\end{problem}

We now proceed to show that this noisy observation SDE insider control problem can be transformed into a \emph{full observation SPDE insider control problem} of the type discussed in the previous sections:\\
To this end, define the probability measure $\tilde{P}$ by
\begin{equation}
d\tilde{P}(\omega)=M_t(\omega) dP(\omega) \text{ on } \mathcal{F}^{v}_t \vee \mathcal{F}_t^{w} \vee \sigma(Z),
\end{equation}
where
\begin{equation}
M_t(\omega)=M_t(\omega,Z)=\exp \Big( -\int_0^t h(X(s,Z)) dw(s) -\frac{1}{2} \int_0^t h^2(X(s,Z))ds \Big).
\end{equation}
It follows by \eqref{Novikov} and the Girsanov theorem that the observation process $R(t)$ defined by \eqref{eq7.2} is a Brownian motion with respect to $\tilde{P}$. Moreover, we have
\begin{equation}
dP(\omega)= K_t(\omega) d\tilde{P}(\omega),
\end{equation}
where
\begin{align}
K_t = M_t^{-1}&= \exp \Big( \int_0^t h(X(s,Z)) dw(s) +\frac{1}{2} \int_0^t h^2(X(s,Z))ds \Big)\nonumber\\
&=\exp \Big( \int_0^t h(X(s,Z)) dR(s) -\frac{1}{2} \int_0^t h^2(X(s,Z))ds \Big).
\end{align}
For $\varphi \in C_0^2(\mathbb{R})$ and fixed $r \in \mathbb{R}, m \in \mathbb{U}$ define the integro-differential operator $L=L_{r,m}$ by
\begin{align}
L_{r,m}\varphi(x)&= \alpha(x,r,m)\frac{\partial \varphi}{\partial x}(x) + \frac{1}{2} \beta^2(x,r,m)\frac{\partial^{2} \varphi}{\partial x^2}\nonumber\\
&+\int_{\mathbb{R}} \{ \varphi(x+\gamma(x,r,m,\zeta))-\varphi(x)-\nabla\varphi(x) \gamma(x,r,m,\zeta) \} \nu(d\zeta),
\end{align}
and let $L^*$ be the adjoint of $L$, in the sense that
\begin{equation}
(L\varphi,\psi)_{L^2(\mathbb{R})}= (\varphi, L^*\psi)_{L^2(\mathbb{R})}
\end{equation}
for all $\varphi,\psi \in C_0^2(\mathbb{R})$.\\

Suppose that for all $z \in \mathbb{R}$ there exists a stochastic process $y(t,x)=y(t,x,z)$ such that
\begin{equation}\label{eq5.8}
\mathbb{E}_{\tilde{P}}[ \varphi(X(t,z)) K_t(z) | \mathcal{R}_t]= \int_{\mathbb{R}} \varphi(x)y(t,x,z)dx
\end{equation}
for all bounded measurable functions $\varphi$.
Then $y(t,x)$ is called the \emph{unnormalized} conditional density of $X(t.z)$ given the observation filtration $\mathcal{R}_t$. Note that by the Bayes rule we have
\begin{equation}\label{eq5.8a}
\mathbb{E}[\varphi(X(t)) | \mathcal{R}_t] = \frac{\mathbb{E}_{\tilde{P}}[\varphi(X(t)) K_t | \mathcal{R}_t]}{\mathbb{E}_{\tilde{P}}[K_t | \mathcal{R}_t]}.
\end{equation}

It is known that under certain conditions the process $y(t,x)=y(t,x,z)$ exists and satisfies the following integro-SPDE, called the Duncan-Mortensen-Zakai equation:

\begin{align}\label{eq7.15}
dy(t,x,z) &= L^*_{R(t),u(t)}y(t,x,z) dt + h(x) y(t,x,z) dR(t); \quad t\geq 0 \nonumber\\
y(0,x,z)&= F(x,z).
\end{align}
See for example Theorem 7.17 in \cite{BC}.

If \eqref{eq5.8} holds, we get
\begin{align}
J(u) &= \mathbb{E}[\int_0^T f(X(t,Z),u(t,Z))dt + g(X(T,Z))]\nonumber\\
&= \mathbb{E}_{\tilde{P}}[\int_0^T f(X(t,Z),u(t,Z))K_t(Z)dt + g(X(T,Z))K_T(Z)] \nonumber\\
&= \mathbb{E}_{\tilde{P}} \big[ \int_0^T \mathbb{E}_{\tilde{P}}[f(X(t,Z),u(t,Z))K_t(Z) | \mathcal{H}_t]dt + \mathbb{E}_{\tilde{P}}[g(X(T,Z)K_T(Z)|\mathcal{H}_T]\Big] \nonumber\\
&= \mathbb{E}_{\tilde{P}} \Big[ \int_0^T \mathbb{E}_{\tilde{P}}[f(X(t,Z),u(t,Z))K_t(Z) | \mathcal{R}_t \vee \sigma(Z)]dt + \mathbb{E}_{\tilde{P}}[g(X(T,Z)K_T(Z) | \mathcal{R}_T \vee \sigma(Z)] \Big]\nonumber\\
&= \mathbb{E}_{\tilde{P}} \Big[ \int_0^T \mathbb{E}_{\tilde{P}}[f(X(t,z),v(z))K_t(z) | \mathcal{R}_t]_{z=Z,v(z)=u(t,z)}dt + \mathbb{E}_{\tilde{P}}[g(X(T,z)K_T(z) |\mathcal{R}_T] _{z=Z} \Big] \nonumber\\
&= \mathbb{E}_{\tilde{P}}[\int_0^T  \int_{\mathbb{R}}f(x,u(t,Z))y(t,x,Z)dx dt + \int_{\mathbb{R}}g(x,Z)y(T,x,Z)dx \Big]=:J_{\tilde{P}}(u).\label{eq5.16}
\end{align}
This transforms the insider partial observation SDE control problem 7.1 into an insider full observation SPDE control problem of the type we have discussed in the previous sections.\\

We summarise what we have proved as follows:

\begin{theorem}[From noisy obs. SDE control to full info. SPDE control]
Assume that \eqref{eq5.8} and \eqref{eq7.15} hold. Then the solution $u^*(t,Z)$ of the noisy observation insider SDE control problem 7.1, consisting of \eqref{eq7.1},\eqref{eq7.2},\eqref{eq7.6},
coincides with the solution $u^*$ of the following (full information) insider SPDE control problem:\\
\begin{problem}
Find $u^{\ast}\in \mathcal{A}$ such that
\begin{equation} \label{eq7.6a}
sup_{u \in \mathcal{A}}J_{\tilde{P}}(u) = J_{\tilde{P}}(u^{\ast}),
\end{equation}
where
\begin{equation}
J_{\tilde{P}}(u) = \mathbb{E}_{\tilde{P}}[\int_0^T  \int_{\mathbb{R}}f(x,u(t,Z))y(t,x,Z)dx dt + \int_{\mathbb{R}}g(x,Z)y(T,x,Z)dx \Big],
\end{equation}
and $y(t,x,Z)$ solves the SPDE
\begin{align}\label{eq7.19}
dy(t,x,Z) &= L^*_{R(t),u(t,Z)}y(t,x,Z) dt + h(x) y(t,x,Z) dR(t); \quad t\geq 0 \nonumber\\
y(0,x,Z)&= F(x,Z).
\end{align}

\end{problem}

\end{theorem}

\section {Examples}
\subsection{Example: Optimal control of a second order SPDE, with control not depending on $x$.}
Consider the following controlled stochastic reaction-diffusion equation:
\begin{align}\label{Wealth}
\begin{cases}
dY(t,x,z)=dY^{\pi}(t,x,z)=[\frac{1}{2}{\frac{\partial^2}{\partial x^2}Y(t,x,z)+\pi(t,z)Y(t,x,z)a_0(t,z)}]dt+\pi(t,z)Y(t,x,z)
b_0(t,z)dB(t)
; \quad t\in[0, T]\\
Y(0,x,z)=\alpha(x) > 0; x \in D,
\end{cases}
\end{align}
with performance functional given by
\begin{equation}
J(\pi):= \mathbb{E}[\int_{D}U(x,Y^{\pi}(T,x,Z),Z)dx]=\int_{\mathbb{R}} j(\pi)dz;
\end{equation}
where
\begin{equation}\label{eq8.2}
j(\pi)=j(\pi,z)= \mathbb{E}[\int_D U(x,Y(T,x,z),z)dx \mathbb{E}[\delta_Z(z)|\mathcal{F}_T]],
\end{equation}
and
 \begin{align}
U(x,y,z)=U(x,y,z,\omega): D\times (0,\infty)\times\mathbb{R}\times \Omega \rightarrow \mathbb{R}\nonumber
\end{align}
is a given utility function, assumed to be concave and $\mathcal{C}^1$ with respect to $y$ and $\mathcal{F}_T$-measurable for each $x,y,z$. Let $\mathcal{A}_{\mathbb{H}}$ be the set of $\mathbb{H}$-adapted controls $\pi(t)$ not depending on $x$ and such that
$$ \mathbb{E} [\int_0^T \pi(t)^2 dt] < \infty.$$
Then it is well-known that the corresponding solution $Y^{\pi}(t,x)$ of \eqref{Wealth} is positive for all $t,x$. See e.g. \cite{Be}.
We study the following problem:
\begin{problem}
Find $\hat{\pi} \in \mathcal{A}_{\mathbb{H}}$ such that
\begin{equation}\label{eq8.3}
\sup_{\pi \in \mathcal{A}_{\mathbb{H}}} j(\pi)= j(\hat{\pi}).
\end{equation}
\end{problem}
This is a problem of the type investigated in the previous sections, in the special case with no jumps and with controls $\pi(t,z)$ not depending on $x$, and we can apply the results there to solve it.

The Hamiltonian \eqref{eq4.1} gets the form, with $u=\pi$,
\begin{equation}\label{eq19}
H(t,x,y,\varphi, \pi,p,q,z)= [\frac{1}{2}\frac{\partial^2}{\partial x^2}\varphi(x) +\pi y a_0(t,z)]p + \pi y b_0(t,z)q,
\end{equation}
while the BSDE \eqref{eq4.2a} for the adjoint processes becomes, keeping in mind that
$$(\frac{\partial^2}{\partial x^2})^*=\frac{\partial^2}{\partial x^2},$$
\begin{align}\label{eq20}
\begin{cases}
dp(t,x,z) &= -[\frac{1}{2}\frac{\partial^2}{\partial x^2}p(t,x,z)+\pi(t,z)\{a_0(t,z)p(t,x,z) + b_0(t,z)q(t,x,z)\}]dt
+ q(t,x,z) dB(t); \quad t\in[0, T]\\
p(T,x,z)&= \frac{\partial U}{\partial y}(x,Y(T,x,z),z)\mathbb{E}[\delta_Z(z)|\mathcal{F}_T]\\
p(t,x,z)&=0 ;  \forall (t,x,z)\in[0,T]\times \partial D\times \mathbb{R}.
\end{cases}
\end{align}

The map
\begin{equation}
 \pi\mapsto \int_DH(t,x,Y(t,x,z),Y(t,.,z),\pi, p(t,x,z),q(t,x,z))dx
 \end{equation}
 is maximal when
 \begin{equation}
 \int_DY(t,x,z)[a_0(t,z)p(t,x,z)+b_0(t,z)q(t,x,z)]dx=0.
\end{equation}
From this we get
\begin{equation}
\int_DY(t,x,z)q(t,x,z)dx=-\frac{a_0(t,z)}{b_0(t,z)}\int_DY(t,x,z)p(t,x,z)dx.
\end{equation}

Let \begin{equation}
\tilde{p}(t,z)=\int_{D}Y(t,x,z)p(t,x,z)dx; \quad t \in [0,T].
\end{equation}
Applying the It\^{o} formula to $Y(t,x,z)p(t,x,z)$ we get:
\begin{align}
&d(Y(t,x,z)p(t,x,z))=p(t,x,z)[(\frac{1}{2}\frac{\partial^2 Y}{\partial x^2}(t,x,z)\nonumber\\
&+\pi(t,z)Y(t,x,z)a_0(t,z))dt+\pi(t,z)Y(t,x,z)b_0(t,z)dB(t)]\nonumber\\
&+Y(t,x,z)[(-\frac{1}{2}\frac{\partial^2 p}{\partial x^2}(t,x,z)-\pi(t,z)(a_0(t,z)p(t,x,z)+b_0(t,z)q(t,x,z)))dt+q(t,x,z)dB(t)]\nonumber\\
&+q(t,x,z)\pi(t,z)Y(t,x,z)b_0(t,z)dt\nonumber\\
&=[\frac{1}{2}\frac{\partial^2 Y}{\partial x^2}(t,x,z)p(t,x,z) -\frac{1}{2}\frac{\partial^2 p}{\partial x^2}(t,x,z)Y(t,x,z)]dt\nonumber\\
&+[\pi(t,z)Y(t,x,z)p(t,x,z)b_0(t,z)+Y(t,x,z)q(t,x,z)]dB(t)
\end{align}
Then  get that the dynamics of $\tilde{p}(t,z)$ is given by
\begin{align}\label{eq25'}
 \begin{cases}
d\tilde{p}(t,z) &=\tilde{p}(t,z)[\pi(t,z)b_0(t,z)-\frac{a_0(t,z)}{b_0(t,z)}]dB(t)\\
\tilde{p}(T,z)&= \int_D\frac{\partial U}{\partial y}(x,Y(T,x,z),z)Y(T,x,z)dx \mathbb{E}[\delta_Z(z)|\mathcal{F}_T].
\end{cases}
\end{align}
Thus we obtain that
\begin{align}\label{eq26}
&\tilde{p}(t,z)=\tilde{p}(0,z) \exp(\int_0^t (b_0(s,z)\pi(s,z)-\frac{a_0(s,z)}{b_0(s,z)})dB(s)\nonumber\\
&-\frac{1}{2} \int_0^t (b_0(s,z)\pi(s,z)-\frac{a_0(s,z)}{b_0(s,z)})^{2}ds),
\end{align}

\noindent for some, not yet determined, constant $\tilde{p}(0,z)$.
In particular, for $t=T$ we get, using \eqref{eq25'},
\begin{align}\label{eq26b}
&\int_D\frac{\partial U}{\partial y}(x,Y(T,x,z),z)Y(T,x,z)dx\mathbb{E}[\delta_Z(z)|\mathcal{F}_T]=\tilde{p}(0,z) \exp(\int_0^T (b_0(t,z)\pi(s,z)-\frac{a_0(s,z)}{b_0(s,z)})dB(s)\nonumber\\
 &-\frac{1}{2} \int_0^T (b_0(s,z)\pi(s,z)-\frac{a_0(s,z)}{b_0(s,z)})^{2}ds).
\end{align}
Now assume that
\begin{equation}\label{eq8.13a}
U(x,y,z)=k(x,z)\ln(y),
\end{equation}
for a given bounded positive $\mathcal{F}_T$-measurable random variable $k(x,z)=k(x,z,\omega)$ such that $K(z):=\int_D k(x,z)dx<\infty$ a.s. for all $z$.
Then equation \eqref{eq26b} becomes
\begin{align}\label{eq8.16}
&K(z)\mathbb{E}[\delta_Z(z)|\mathcal{F}_T]=\tilde{p}(0,z) \exp(\int_0^T (b_0(s,z)\pi(s,z)-\frac{a_0(s,z)}{b_0(s,z)})dB(s)\nonumber\\
 &-\frac{1}{2} \int_0^T (b_0(t,z)\pi(s,z)-\frac{a_0(s,z)}{b_0(s,z)})^{2}ds).
\end{align}

\noindent To make this more explicit, we proceed as follows:\\

Define
\begin{equation}
M(t,z) := \mathbb{E}[K(z)\mathbb{E}[\delta_Z(z)|\mathcal{F}_T] | \mathcal{F}_t];\quad 0 \leq t \leq T.
\end{equation}
Then by the generalized Clark-Ocone theorem in \cite{AaOPU},
\begin{equation}
\begin{cases}
dM(t,z) = \mathbb{E}[D_t M(T,z) | \mathcal{F}_t]dB(t)= \Phi_K(t,z)M(t,z) dB(t)\\
M(0,z)=1
\end{cases}
\end{equation}
where
\begin{equation}\label{eq18a}
\Phi_K(t,z)= \frac{\mathbb{E}[D_tM(T,z)|\mathcal{F}_t]}{M(t,z)}=\frac{\mathbb{E}[D_t[K(z)\mathbb{E}[\delta_Z(z)|\mathcal{F}_T]] | \mathcal{F}_t]}{\mathbb{E}[K(z)\mathbb{E}[\delta_Z(z)|\mathcal{F}_T] | \mathcal{F}_t]}.
\end{equation}
Solving this SDE for $M(t,z)$ we get
\begin{equation}
M(t,z) = \exp( \int_0^t \Phi_K(s,z)dB(s) -\frac{1}{2} \int_0^t \Phi_K^2(s,z) ds).
\end{equation}
Since the two martingales $\tilde{p}(t,z)$ and $M(t,z)$ are identical for $t=T,$ they are identical for all $t \leq T$ and hence, by \eqref{eq8.16} we get
\begin{align}
&M(0,z)\exp( \int_0^t \Phi_K(s,z)dB(s) -\frac{1}{2} \int_0^t \Phi_K^2(s,z) ds)\nonumber\\
&=\tilde{p}(0,z) \exp(\int_0^t (b_0(s,z)\pi(s,z)-\frac{a_0(s,z)}{b_0(s,z)})dB(s)\nonumber\\
&-\frac{1}{2} \int_0^t (b_0(s,z)\pi(s,z)-\frac{a_0(s,z)}{b_0(s,z)})^{2}ds);\quad 0\leq t \leq T.
\end{align}
By identification of the integrals with respect to $ds$ and the stochastic integrals we get:
\begin{equation}\label{eq8.21}
\pi(t,z)=\hat{\pi}(t,z)=\frac{\Phi_K(t,z)}{b_0(t,z)}+\frac{a_0(t,z)}{\sigma^2_0(t,z)}.
\end{equation}
We summarise what we have proved as follows:
\begin{theorem}
The optimal insider control $\hat{\pi}$ for the problem \eqref{eq8.3} with $U(x,y,z)=k(x,z)\ln(y)$ as in \eqref{eq8.13a}, is given by \eqref{eq8.21}, with $\Phi_K(t,z)$ given by \eqref{eq18a}.
\end{theorem}

\begin{coro}
Suppose $k(x,z)$ is deterministic. Then the optimal insider control $\hat{\pi}$ for the problem \eqref{eq8.3} with $U(x,y,z)=k(x,z)\ln(y),$ is given by
\begin{equation}
\hat{\pi}(t,z)=\hat{\pi}_1(t,z)=\frac{\Phi_1(t,z)}{b_0(t,z)}+\frac{a_0(t,z)}{\sigma^2_0(t,z)},
\end{equation}
with
\begin{equation}
\Phi_1(t,z)=\frac{\mathbb{E}[D_t\delta_Z(z)|\mathcal{F}_t]}{\mathbb{E}[\delta_Z(z)|\mathcal{F}_t]}.
\end{equation}
Note that the optimal insider portfolio in this case is in fact the same as in the case when $Y$ does not depend on $x$. See \cite{DO1}.
\end{coro}

\subsection{Optimal insider portfolio with noisy observations}

We now study an example illustrating the application in Section 7: \\
Let $\alpha$ and $\beta$ be given adapted processes, with $\beta$ bounded away from 0. Suppose the signal process $X(t)=X^{\pi}(t,Z)$ is given by
\begin{equation} \label{eq8.18}
\begin{cases}
dX(t,Z) = \pi(t,Z)[\alpha(t) dt + \beta(t) dv(t)]; \quad 0 \leq t \leq T,\\
X(0) \text{ has density } F
\end{cases}
\end{equation}
Here $\pi(t,Z)$ is the control, representing the portfolio in terms of the amount invested in the risky asset at time $t$, when the risky asset unit price $S(t)$ is given by
\begin{equation}
\begin{cases}
dS(t)= S(t)[\alpha(t) dt + \beta(t) dv(t)]; \quad 0 \leq t \leq T,\\
S(0) > 0,
\end{cases}
\end{equation}
and the safe investment unit price is $S_0(t)=1$ for all $t$. The process $X(t)$ then represents the corresponding value of the investment at time $t$. For $\pi$ to be in the set $\mathcal{A}_{\mathbb{H}}$ of admissible controls, we require that $X(t) > 0$ for all $t$ and
\begin{equation}
\mathbb{E}[\exp(\frac{1}{2}\int_0^TX^2(s,Z)ds)] < \infty.
\end{equation}
See \eqref{Novikov}.
Suppose the observations $R(t)$ of $X(t)$ at time $t$ are not exact, but subject to uncertainty or noise, so that
the observation process is given by
\begin{equation}
\begin{cases}
dR(t,Z)=  X(t,Z)dt + dw(t); \quad 0 \leq t \leq T,\\
R(0)=0.
\end{cases}
\end{equation}
Here, as in Section 7, the processes $v$ and $w$ are independent Brownian motions, and the random variable $Z$ represents the information available to the insider from time 0.
Let $U: [0,\infty) \mapsto [-\infty,\infty)$ be a given $C^{1}$ (concave) utility function. The performance functional is assumed to be
\begin{equation}
J(\pi)= \mathbb{E}[U(X^{\pi}(T,Z))].
\end{equation}
By Theorem 7.2, the problem to maximize $J(\pi)$ over all $\pi \in \mathcal{A}_{\mathbb{H}}$  is equivalent to the following problem:\\

\begin{problem}
Find $\hat{\pi} \in \mathcal{A}$ such that
\begin{equation}\label{eq8.23}
sup_{\pi \in \mathcal{A}}J_{\tilde{P}}(\pi)=J_{\tilde{P}}(\hat{\pi}),
\end{equation}
where
\begin{equation}
J_{\tilde{P}}(\pi)=\mathbb{E}_{\tilde{P}}[\int_{\mathbb{R_{+}}}U(x,Z)y(T,x,Z)dx],
\end{equation}
and $y(t,x,Z)=y^{\pi}(t,x,Z)$ is the solution of the SPDE
\begin{equation}\label{eq8.24}
\begin{cases}
dy(t,x,Z) = (L^*_{\pi(t,z)}y)(t,x,Z)dt +x y(t,x,Z) dR(t); \quad 0 \leq t \leq T,\\
y(0,x,Z)=F(x),
\end{cases}
\end{equation}
\end{problem}
where
\begin{equation}\nonumber
(L^*_{\pi(t)}y)(t,x,z)=-\pi(t,z)\alpha(t) y'(t,x,z) + \frac{1}{2} \pi^2(t,z)\beta^2(t) y''(t,x,z),
\end{equation}
with $y'(t,x,z)=\frac{\partial y(t,x,z)}{\partial x}, y''(t,x,z)=\frac{\partial^2 y(t,x,z)}{\partial x^2}.$\\

Define the space
\begin{equation}
\mathbf{H}^1(\mathbb{R}^+)=\{y\in\mathbf{L}^2(\mathbb{R}^+), \frac{\partial y}{\partial x}\in\mathbf{L}^2(\mathbb{R}^+)\}
\end{equation}
The $\mathbf{H}^1$ norm is given by:
\begin{equation}
\|y(t,z)\|^2_{\mathbf{H}^1 (\mathbb{R}^+)}=\|y(t,z)\|^2_{\mathbf{L}^2(\mathbb{R}^+)}+\|y'(t,z)\|^2_{\mathbf{L}^2(\mathbb{R}^+)}
\end{equation}
We have
\begin{equation}
\mathbf{H}^1(\mathbb{R}^+)\subset\mathbf{L}^2(\mathbb{R}^+)\subset \mathbf{H}^{-1}(\mathbb{R}^+)
\end{equation}
We verify the coercivity condition of the operator $-L^*_{\pi(t)}$:
\begin{align}
2\langle -L^*_{\pi(t)}y, y\rangle&=2\pi(t,z)\alpha(t) \langle y'(t,x,z), y(t,x,z)\rangle - \pi^2(t,z)\beta^2(t) \langle y''(t,x,z)y(t,x,z)\rangle\nonumber\\
&=2\pi(t,z)\alpha(t) \int_{\mathbb{R}^+} y'(t,x,z) y(t,x,z)dx- \pi^2(t,z)\beta^2(t)\int_{\mathbb{R}^+} y''(t,x,z)y(t,x,z)dx\nonumber\\
&=\pi(t,z)\alpha(t)[y^2(t,x,z)]_{\partial \mathbb{R}^+}-\pi^2(t,z)\beta^2(t)[y(t,x,z)y'(t,x,z)]_{\partial \mathbb{R}^+}\nonumber\\
&+\pi^2(t,z)\beta^2(t)\int_{\mathbb{R}^+} (y'(t,x,z))^2dx.
\end{align}
Suppose that $y(t,x,z)=0$ for $x=0$. Then we get
\begin{equation}
2\langle -L^*_{\pi(t)}y, y\rangle=\pi^2(t,z)\beta^2(t)\|y'(t,z)\|^2_{\mathbf{L}^2(\mathbb{R}^+)}.
\end{equation}
Let
\begin{equation}
\mathbf{H}^1_0(\mathbb{R}^+)=\{y\in\mathbf{H}^1, y=0 \text{ on } \partial \mathbb{R}^+\}.
\end{equation}

We have $|y(t,z)|_{1,\mathbb{R}^+}=\|y'(t,z)\|_{\mathbf{L}^2(\mathbb{R}^+)}$ is a norm in  $\mathbf{H}^1_0(\mathbb{R}^+)$, which is equivalent to the $\mathbf{H}^1 (\mathbb{R}^+)$ norm;  i.e. there exist $a,b>0$ such that
\begin{equation}
a\|y(t,z)\|_{1,\mathbb{R}^+}\leq |y(t,z)|_{1,\mathbb{R}^+}=\|y'(t,z)\|_{\mathbf{L}^2(\mathbb{R}^+)}\leq b\|y(t,z)\|_{1,\mathbb{R}^+}
\end{equation}
We conclude that the following coercivity condition is satisfied:
\begin{equation}
2\langle -L^*_{\pi(t)}y, y\rangle\geq a^2\pi^2(t,z)\beta^2(t)\|y(t,z)\|_{1,\mathbb{R}^+}^2.
\end{equation}
Using Theorem 1.1 and Theorem 2.1 in Pardoux \cite{Par}, we obtain that \eqref{eq8.24} has a unique solution $y(.,.,z)\in \mathbf{L}^2(\Omega,\mathbf{C}(0,T,\mathbf{L}^2(\mathbb{R}^+)))$ i.e. $y(.,.,z)$ satisfies
\begin{enumerate}
\item  $\mathbb{E}[y^2(t,x,z)] < \infty$ for all $t,x,z$.
\item  The map $t \mapsto y(t,.,z)$ is continuous as a map from $[0,T]$ into $L^2(\mathbb{R}^{+})$, for all $z$.
\end{enumerate}

 Moreover, the first and second partial derivatives with respect to $x$, denoted by $y'(t,x,z)$ and $y''(t,x,z)$ respectively, exist and belong to $L^2(\mathbb{R})$.

The problem \eqref{eq8.23} is of the type discussed in Section 6 and we now apply the methods developed there to study it:\\
The Hamiltonian given in \eqref{eq4.1} now gets the form
\begin{equation}
H(t,x,y,\varphi,\pi,p,q)=(L^*_{\pi}\varphi) p + x y q,
\end{equation}
and the adjoint BSDE \eqref{eq4.2a} becomes
\begin{equation}\label{eq8.37}
\begin{cases}
dp(t,x,z)= -[A_{\pi(t,z)}p(t,x,z)+xq(t,x,z)] dt + q(t,x,z) dR(t); \quad 0 \leq t \leq T,\\
p(T,x,z)=U(x,z)\mathbb{E}_{Q}[\delta_{Z}(z)|\mathcal{R}_T].
\end{cases}
\end{equation}
where $\mathcal{R}_t$ is the sigma-algebra generated by $\{R(s)\}_{s\leq t}$, for $0 \leq t \leq T,$ and
\begin{equation}\label{eq8.43}
A_{\pi(t,z)}p(t,x,z)=\pi(t,z)\alpha(t) p'(t,x,z) + \frac{1}{2} \pi^2(t,z)\beta^2(t) p''(t,x,z).
\end{equation}

By \cite{OPZ} and \cite{ZRW}, this backward SPDE (BSPDE for short) admits a unique solution which belongs to $\mathbf{L}^2(\mathbb{R}^+)$.


The map
$$ \pi \mapsto \int_{\mathbb{R}^{+}}H(t,x,y(t,x,z),y(t,.,z),\pi,p(t,x,z),q(t,x,z))dx$$
 is maximal when
\begin{equation}
\int_{\mathbb{R}^{+}} \{- \alpha(t) y'(t,x,z)+\pi \beta^2(t)y''(t,x,z) \}p(t,x,z)dx =0,
\end{equation}
i.e. when $\pi=\hat{\pi}(t,z)$, given by
\begin{equation} \label{eq8.29}
\hat{\pi}(t,z)=\frac{\alpha(t)\int_{\mathbb{R}_+}y'(t,x,z)p(t,x,z)dx}{\beta^2(t) \int_{\mathbb{R}_+} y''(t,x,z)p(t,x,z)dx}.
\end{equation}
Using integration by parts and \eqref{eq5.8} -\eqref{eq5.8a} we can rewrite this as follows:
\begin{align}\label{eq8.29a}
\hat{\pi}(t,z)&=-\frac{\alpha(t)\int_{\mathbb{R}_+}y(t,x,z)p'(t,x,z)dx}{\beta^2(t)\int_{\mathbb{R}_+} y(t,x,z)p''(t,x,z)dx}\nonumber\\
&=-\frac{\alpha(t)\mathbb{E}[p'(t,X(t),z)|\mathcal{R}_t]}{\beta^2(t)\mathbb{E}[p''(t,X(t),z)|\mathcal{R}_t]}.
\end{align}
We summarise what we have proved as follows:
\begin{theorem}
Assume that the conditions of Theorem 7.2 hold. A portfolio $\hat{\pi}(t,z) \in \mathcal{A}$ is an optimal portfolio for the noisy observation insider portfolio problem \eqref{eq8.23}, if it is given in feedback form by 
\begin{equation}\label{eq8.46}
\hat{\pi}(t,z)=-\frac{\alpha(t)\mathbb{E}[p'(t,X(t),z)|\mathcal{R}_t]}{\beta^2(t)\mathbb{E}[p''(t,X(t),z)|\mathcal{R}_t]},
\end{equation}
 where $p(t,x,z)$ solves the BSPDE
\begin{equation}\label{eq8.47}
\begin{cases}
dp(t,x,z)= -[A_{\hat{\pi}(t,z)}p(t,x,z)+xq(t,x,z)] dt + q(t,x,z) dR(t);0\leq t \leq T,\\
p(T,x,z)=U(x)\mathbb{E}_{\tilde{P}}[\delta_{Z}(z)|\mathcal{R}_T],
\end{cases}
\end{equation}
and $p''(t,x,z) \neq 0$ for all $t,x,z$.
\end{theorem}


\begin{thebibliography}{9999}

\bibitem[Aa\O PU]{AaOPU}
K. Aase, B. \O{}ksendal, N. Privault and J. Ub\o e: White noise generalizations of the Clark-Haussmann-Ocone theorem with application to mathematical finance. Finance Stoch. 4 (2000), 465-496.

\bibitem[Aa\O U]{AaOU}
K. Aase, B. \O{}ksendal and J. Ub\o e: Using the Donsker delta function to compute hedging strategies. Potential Analysis 14 (2001), 351-374.

\bibitem[BC]{BC} A. Bain and D. Crisan: Fundamentals of Stochastic Filtering. Springer 2009.

\bibitem[Be]{Be} F.E. Benth: On the positivity of the stochastic heat equation. Potential Analysis 6 (1997), 127-148.

\bibitem[B]{B} L. Breiman: Probability. Addison-Wesley 1968.

\bibitem[B\O]{BO}
F. Biagini and B. \O{}ksendal: A general stochastic calculus approach to insider trading.Appl. Math. \& Optim. 52 (2005), 167-181.

\bibitem[DM\O P1]{DMOP1} G. Di Nunno, T. Meyer-Brandis, B. {\O}ksendal and F. Proske: Malliavin calculus and anticipative It\^ o formulae for L\'{e}vy processes. Inf. Dim. Anal. Quantum Prob. Rel. Topics 8 (2005), 235-258.

\bibitem[DM\O P2]{DMOP2} G. Di Nunno, T. Meyer-Brandis, B. {\O}ksendal and F. Proske: Optimal portfolio for an insider in a market driven by L\'{e}vy processes. Quant. Finance 6 (2006), 83-94.

\bibitem[D\O 1]{DO1} O. Draouil and B. \O ksendal: A Donsker delta functional approach to optimal insider control and application to finance. Comm. Math. Stat.  (CIMS) 3 (2015), 365-421; DOI 10.1007/s40304-015-0065-y.

\bibitem[D\O 2]{DO3} O. Draouil and B. \O ksendal:  Optimal insider control and semimartingale decompositions under enlargement of filtration. arXiv: 1512.01759v1 (6 Dec.2015). To appear in Stochastic Analysis and Applications.

\bibitem[Di\O 1]{DiO1} G. Di Nunno and B. {\O}ksendal: The Donsker delta function, a representation formula for functionals of a L\'{e}vy process and application to hedging in incomplete markets. S\'{e}minaires et Congr\`{e}es, Societ\'{e} Math\'{e}matique de France, Vol. 16 (2007), 71-82.

\bibitem[Di\O 2]{DiO2} G. Di Nunno and B. \O ksendal: A representation theorem and a sensitivity result for functionals of jump diffusions. In A.B. Cruzeiro, H. Ouerdiane and N. Obata (editors): Mathematical Analysis and Random Phenomena. World Scientific 2007, pp. 177 - 190.

\bibitem[D\O P]{DOP} G. Di Nunno, B. {\O}ksendal and F. Proske: Malliavin Calculus for L\'{e}vy Processes with Applications to Finance. Universitext, Springer 2009.

\bibitem[H\O UZ]{HOUZ} H. Holden, B. {\O}ksendal, J. Ub\o e and T. Zhang: Stochastic Partial Differential Equations. Universitext, Springer, Second Edition 2010.

\bibitem[LP]{LP} A. Lanconelli and F. Proske: On explicit strong solution of It\^ o-SDEs and the Donsker delta function of a diffusion.
Inf. Dim. Anal. Quatum Prob Rel. Topics 7 (2004),437-447.

\bibitem[M\O P]{MOP} S. Mataramvura, B. \O ksendal and F. Proske: The Donsker delta function of a L\'{e}vy process with application to chaos expansion of local time. Ann. Inst H. Poincar\'{e} Prob. Statist. 40 (2004), 553-567.

\bibitem[MP]{MP} T. Meyer-Brandis and F. Proske: On the existence and explicit representability of strong solutions of L\'{e}vy noise driven SDEs with irregular coefficients. Commun. Math. Sci. 4 (2006), 129-154.

\bibitem[\O1]{O1} B. \O ksendal: Optimal control of stochastic partial differential equations.
Stochastic Analysis and Applications 23 (2005), 165-179.

\bibitem[\O PZ]{OPZ} B. \O ksendal, F. Proske and T. Zhang:
Backward stochastic partial differential equations with jumps and application to optimal control of random jump fields.
Stochastics 77 (2005), 381-399.


\bibitem[\O S1]{OS1}
B. \O{}ksendal and A. Sulem: Applied Stochastic Control of Jump Diffusions. Second Edition. Springer 2007

\bibitem[\O S2]{OS2}
B. \O{}ksendal and A. Sulem: Risk minimization in financial markets modeled by It\^ o-L\' evy processes. Afrika Matematika (2014), DOI: 10.1007/s13370-014-02489-9.

\bibitem[Par]{Par} E. Pardoux: Stochastic Partial Differential Equations and filtering of diffusion processes. Stochastics. 1979, Vol3,pp-127-167.

\bibitem[P]{P} P. Protter: Stochastic Integration and Differential Equations. Second Edition. Springer 2005

\bibitem[PK]{PK} I. Pikovsky and I. Karatzas: Anticipative portfolio optimization. Adv. Appl. Probab. 28 (1996), 1095-1122.

\bibitem[PR]{PR} CI. Pr\'{e}v\^{o}t and M. Roeckner: A concise course on stochastic partial differential equations. Lecture Notes in Mathematics 1905, Springer 2007.

\bibitem[RV]{RV} F. Russo and P. Vallois: Forward, backward and symmetric stochastic integration. Probab. Theor. Rel. Fields 93 (1993), 403-421.

\bibitem[RV1]{RV1} F. Russo and P. Vallois. The generalized covariation process and It\^{o} formula. Stoch. Proc. Appl., 59(4):81-104, 1995.

\bibitem[RV2]{RV2} F. Russo and P. Vallois. Stochastic calculus with respect to continuous finite quadratic variation processes. Stoch. Stoch. Rep., 70(4):1-40, 2000.

\bibitem[ZRW]{ZRW} Q. Zhou, Y. Ren and W. Wu. On solutions to backward stochastic partial differential equations for L\'{e}vy processes. Journal of Computational and Applied Mathematics 235 (2011), 5411-5421.

\end{thebibliography}
\end{document}